\documentclass[12pt,reqno]{amsart}

\usepackage{amsfonts}
\usepackage{amssymb}
\usepackage{amsthm}
\usepackage{amsmath}

\newtheorem{theo}{Theorem}[section]

\newtheorem{coro}[theo]{Corollary}
\newtheorem{lem}[theo]{Lemma}

\theoremstyle{definition}
\newtheorem*{defin}{Definition}
\newtheorem*{notat}{Notation}

\theoremstyle{remark}
\newtheorem*{rema}{Remark}

\numberwithin{equation}{section}

\DeclareMathOperator{\md}{md}
\DeclareMathOperator{\mdn}{\overline{md}}

\newcommand{\ga}{>_{\text{a}}}
\newcommand{\gea}{\ge_{\text{a}}}
\newcommand{\la}{<_{\text{a}}}

\newcommand{\ebar}{\overline{e}}
\newcommand{\Mbar}{\overline{M}}

\newcommand{\bs}{\boldsymbol}

\newcommand{\T}{^{\rm T}}
\newcommand{\dmin}{d_{\text{min}}}
\newcommand{\emax}{e_{\text{max}}}
\newcommand{\Ftilde}{{\widetilde F}}
\newcommand{\uvec}{{\bs{\mathrm u}}}
\newcommand{\vvec}{{\bs{\mathrm v}}}

\begin{document}

\title[The maximal $\{-1,1\}$-determinant of order~$15$]{The maximal
$\{-1,1\}$-determinant of order~$15$}
\author[W. P. Orrick]{W. P. Orrick}
\address{Department of Mathematics, Indiana University,
Bloomington IN 47405, USA}
\date{\today}

\begin{abstract}
We study the question of finding the maximal determinant of matrices of
odd order with entries $\pm1$.  The most general upper bound on the maximal
determinant, due to Barba,
can only be achieved when the order is the sum of two consecutive squares.
It is conjectured that the bound is always attained in
such cases.  Apart from these, only in orders 3, 7, 9, 11, 17 and 21 has
the maximal value been established.  In this paper we confirm the results
for these orders, and
add order 15 to the list.  We follow previous authors in
exhaustively searching for candidate Gram matrices having determinant
greater than or equal to the square of a known lower bound on the
maximum.  We then attempt to decompose each candidate as the product of a
$\{-1,1\}$-matrix and its transpose.  For order 15 we find four candidates, all
of Ehlich block form, two having determinant $(105\cdot3^5\cdot2^{14})^2$ and
the others determinant $(108\cdot3^5\cdot2^{14})^2$.  One of the former
decomposes (in an essentially unique way)
while the remaining three do not.  This result proves a
conjecture made independently by W. D. Smith and J. H. E. Cohn.  We also
use our method to compute improved upper bounds on the maximal determinant
in orders 29, 33, and 37, and to establish the range of the determinant
function of $\{-1,1\}$-matrices in orders 9 and 11.
\end{abstract}

\maketitle

\section{Introduction}\label{sect:intro}
Under the conditions that its elements are real and bounded in magnitude,
what is the largest determinant that a square matrix can have?  By rescaling
we may fix the bound on the magnitude of the elements to 1, and because of
linearity of the determinant in its rows and columns, we lose nothing by
restricting
the entries to the set $\{-1,1\}$.  For $R_n\in\{\text{$n\times n$ matrices
with entries $\pm1$}\}$, Hadamard~\cite{had} showed that
\begin{equation}
\det R_n\le n^{n/2}
\label{Bzero}
\end{equation}
and that equality can only hold when $n=1$, $2$ or $n\equiv0\bmod 4$.  It is
not known whether equality always does hold in such cases, but
Paley~\cite{pal} conjectured
that it does, and the lowest order for which the question is unresolved is
$n=428$.  For other orders more stringent bounds have been established.
Barba~\cite{bar} showed that when $n$ is odd the bound is
\begin{equation}
\det R_n\le\sqrt{2n-1}\ (n-1)^{(n-1)/2}:=B(n).
\label{Bone}
\end{equation}
This bound is an integer only when $n$ is the sum of two consecutive squares,
and it is believed that equality holds whenever this is the case.  However,
even for orders as low as $85=6^2+7^2$ this has not been proved.  When
$n\equiv2\bmod4$, the bound, derived independently by Ehlich~\cite{ehla}
and Wojtas~\cite{Wo}, is
\begin{equation}
\det R_n\le 2(n-1)(n-2)^{(n-2)/2}.
\label{Btwo}
\end{equation}
The bound has been shown to be attainable only if $n-1$ is the sum of two
squares~\cite{ehla}, and again it is believed to hold in all such cases.
The lowest order for which the question has not been settled is
$n=138$~\cite{FKS}.

We say that a bound is {\it tight} when it is achieved infinitely often.
All three of the above bounds have been shown constructively to be
tight~\cite{bro,spe,kks}.  The
bound~(\ref{Bone}) is actually tight only for $n\equiv1\bmod4$.  For
$n\equiv3\bmod4$ Ehlich~\cite{ehlb} proved the stricter bound
\begin{equation}
\det R_n\le(n-3)^{(n-s)/2}(n-3+4r)^{u/2}(n+1+4r)^{v/2}
\sqrt{1-\frac{ur}{n-3+4r}-\frac{v(r+1)}{n+1+4r}}
\label{Bthree}
\end{equation}
where $s=5$ for $n=7$, $s=5$ or $6$ for $n=11$, $s=6$ for
$15\le n\le59$ and s=7 for $n\ge63$ and where
$r=\left\lfloor\frac{n}{s}\right\rfloor$,
$n=rs+v$ and $u=s-v$.  This bound is not integral for $n\le59$.  For $n\ge63$
Cohn~\cite{Co2} has shown that it is integral only when $n=112t^2\pm28t+7$ for
some integer $t$, but even in these cases, it is not known whether the bound
is ever attained.

Denote by $\md(n)$ the maximal determinant attained by any member of the
set of  $n\times n$ $\{-1,1\}$-matrices and define $\mdn(n):=\md(n)/2^{n-1}$.
The latter is the ``normalized'' maximal determinant function where the
factor $2^{n-1}$ common to all $n\times n$ $\{-1,1\}$-matrices has been
removed.

When a determinant is found that attains the relevant one of the above
bounds, it is immediate that $\md(n)$ is just the bound
itself, but when the upper bound is not attained, finding $\md(n)$
can be exceedingly difficult.  To date, the only orders in this category for
which $\md(n)$ is known are 3, 7~\cite{Wi}, 9~\cite{EZ}, 11~(Ehlich,
as reported in~\cite{GK}), 17~\cite{MK} and 21~\cite{CKM}.  This leaves, for
$n\le30$, orders 15, 19, 22, 23, 27, and 29 unresolved.  The sequence of
maximal determinants for $1\le n\le14$ is
\begin{equation*}
\begin{array}{r|cccccccccccccc}
n & 1 & 2 & 3 & 4 & 5 & 6 & 7 & 8 & 9 & 10 & 11 & 12 & 13 & 14 \\
\mdn(n) & 1 & 1 & 1 & 2 & 3 & 5 & 9 & 32 & 56 & 144 & 320 & 1458 & 3645 &
9477
\end{array}\,.
\end{equation*}
In this paper, we settle the case $n=15$.  This means
the above sequence can be extended up to $n=18$:
\begin{equation*}
\begin{array}{r|cccc}
n & 15 & 16 & 17 & 18 \\
\mdn(n) & 25515 & 131072 & 327680 & 1114112
\end{array}\,.
\end{equation*}

Our method of proof follows
that of Ehlich (reported in~\cite{GK}), Moyssiadis and Kounias~\cite{MK}
and Chadjipantelis, Kounias and Moyssiadis~\cite{CKM}, and relies on an
intensive computer search.  Compared with all other orders in which $\md(n)$
has been determined, a considerably greater (but still modest) amount of
computer power is required for $n=15$.  The value 25515 was
conjectured to be the maximum by Warren D. Smith~\cite{Sm} and independently
by John H. E. Cohn~\cite{Co1,Co2}.

In addition to this result for $n=15$, we confirm the known values of $\md(n)$
for all odd $n\le21$ excluding, of course, $n=19$.  We are also able to report
improved upper
bounds on the maximal determinants for matrices of orders 29, 33, and 37.
Furthermore, we establish the range of the determinant function for matrices
of orders 9 and 11.

\section{Some properties of matrices}\label{sect:matProps}
While we are here mainly interested in the case $n=15$, we also would like
to confirm existing results in other orders and to say as much as possible
about other open cases.  On the other hand, even orders (the lowest
interesting case is $n=22$) require a somewhat different approach, owing
primarily to the absence of a canonical normalization (see below).
Hence we will describe an algorithm that works for odd orders generally.
For the remainder of this paper we take $n$ to be an odd positive number.

\subsection{Equivalences of $\{-1,1\}$-matrices}
Negating some set of rows and some set of columns of a matrix changes its
determinant by at most a sign.  Two matrices related to each other by such
negations are considered to be equivalent.  It is advantageous
to restrict our attention to
a single canonical representative of each equivalence class:
\begin{defin}
A vector of odd length with elements in the set $\{-1,1\}$ is
{\em parity normalized} if it has an even number
of positive elements.  A matrix of odd order with elements in $\{-1,1\}$ is
parity normalized all its rows and columns are.
\end{defin}
\begin{lem}
Any $\{-1,1\}$-matrix of odd order may be converted to a unique parity
normalized matrix by a series of negations of rows and columns.
\end{lem}
\begin{proof}
Begin by negating all rows containing an even number of entries $-1$.
Since the number of rows is odd, the matrix now contains an
odd number of negative entries.  Thus the number of columns containing
an even number of negative entries must be even.  Negating these gives
the desired matrix, which is clearly unique.
\end{proof}

Permuting a set of rows and a set of columns of a matrix also preserves
its determinant up to a sign. 
\begin{defin}
Two $\{-1,1\}$-matrices $R$ and $S$ are {\em Hadamard equivalent} if
$S=PRQ$ for some pair of signed permutation matrices $(P,Q)$.
\end{defin}
Note that we have omitted a third determinant preserving operation,
transposition, in our notion of equivalence.  Henceforth, we take
$\{-1,1\}$-matrices to be parity normalized.  On the other hand, transforming
a given matrix into a canonical representative of the same class under
full Hadamard equivalence is computationally expensive, so we do not attempt
to impose a canonical ordering on the rows and columns of a matrix.

\subsection{Gram matrices}\label{subsect:GramProperties}
Let $M_r(R)$ be the Gram matrix of the rows of $R$,
\begin{equation}
M_r(R)=RR^{\rm T}.
\end{equation}
$M_r(R)$ is unchanged under negation and permutation of the columns of $R$.
Likewise let $M_c(R)$ be the Gram matrix of the columns of $R$,
\begin{equation}
M_c(R)=R^{\rm T}R.
\end{equation}
$M_c(R)$ is unchanged under negation or permutation of the rows of $R$.
Both $M_r(R)$ and $M_c(R)$ are symmetric, have diagonal elements equal to $n$,
and, assuming $R$ nonsingular, are positive definite.  They have the
same determinant, which must be a perfect square.  Furthermore, they have
identical characteristic equations.  This follows from the following lemma.
\begin{lem}\label{lemma:charEqn}
Let $A$ and $B$ be square matrices.  Then
$\det(AB-\lambda I)=\det(BA-\lambda I)$.
\end{lem}
\begin{proof}
If $B$ is invertible, the result follows immediately by multiplying the
right-hand side of the equation on the left by $\det B$ and on the right
by $\det B^{-1}$, and then combining determinants.  If $B$ is not
invertible, consider $B'=B-x I$.  Since
$B'$ is invertible for all but finitely many $x$, we have, in the generic
case, that $\det(AB'-\lambda I)=\det(B'A-\lambda I)$.  As this
is a polynomial identity in $x$, it must be true for all $x$ and for
$x=0$ in particular.
\end{proof}

Since Gram matrices are symmetric, we will often have occasion to talk about
operations that act on rows and columns simultaneously in the same way.
We will therefore refer to such operations as, for example, permuting the
{\em indices} of $M$.

When $R$ is parity normalized, the elements of $M_r(R)$ and $M_c(R)$ are
necessarily congruent to $n\bmod4$.  Under permutation of the
rows of $R$, the indices of $M_r(R)$ undergo the same permutation.  Under
permutation of the columns of $R$, the indices of $M_c(R)$ undergo the same
permutation.

The search method we employ requires the generation of a set of candidate
matrices which are potentially Gram matrices of $\{-1,1\}$-matrices.  We are
interested in $\{-1,1\}$-matrices whose determinant equals or exceeds some
threshold, $\dmin$ in magnitude.  For a given $n$, this threshold is usually
chosen to be the best known lower bound on $\md(n)$.  To be a candidate, a
matrix $M$ must possess
the properties of $M_r(R)$ and $M_c(R)$ described above and must have
determinant greater than or equal to $\dmin^2$.
The set of such candidates is the set ${\mathcal M}_n^\square(\dmin)$ whose
definition is made precise below.

\begin{defin} We denote by ${\mathcal M}_{n,p}$ be the set of $p\times p$
matrices $M$ having the following properties:
\begin{enumerate}
\item $M$ is symmetric and positive definite;
\item the diagonal elements of $M$ are all equal to $n$;
\item the elements of $M$ are integers congruent to $n \bmod{4}$.
\end{enumerate}
The set ${\mathcal M}_{n,p}(\dmin)$ is the subset of
${\mathcal M}_{n,p}$ whose elements satisfy the additional requirement
\begin{enumerate}
\setcounter{enumi}{3}
\item $\det{M}\ge \dmin^2$.
\end{enumerate}
We define ${\mathcal M}_n={\mathcal M}_{n,n}$ and
${\mathcal M}_n(\dmin)={\mathcal M}_{n,n}(\dmin)$.  Finally, let
${\mathcal M}_n^\square(\dmin)$ be the subset of ${\mathcal M}_n(\dmin)$
whose elements satisfy
\begin{enumerate}
\setcounter{enumi}{4}
\item $\det{M}=\text{perfect square.}$
\end{enumerate}
\end{defin}
Positive definiteness implies that the off-diagonal elements of $M$ are
less than $n$ in absolute value.

Just as for $\{-1,1\}$-matrices, there is a notion of equivalence for
candidate Gram matrices.  Negation will not be considered because
the congruence property modulo 4 that we have imposed on the matrix
elements already fixes the signs.
\begin{defin} Two matrices, $M_1,M_2\in{\mathcal M}_{n,p}$, are
{\em equivalent} if one can be obtained from the other by a
permutation of indices.  The equivalence class of $M$
will be denoted ${\mathcal E}_M$.
\end{defin}
\begin{defin}An element of a matrix $M\in {\mathcal M}_{n,p}$ is
{\em minimal} if its magnitude is the smallest allowed value.  If
$n\equiv1\bmod4$ the minimal element is $1$; if $n\equiv3\bmod4$ the
minimal element is $-1$.  The {\em minimal vector}
of length $k$ is the vector of length $k$, all of whose elements are minimal.
\end{defin}
\begin{defin}A set of indices $A\subseteq \{1,\ldots,p\}$ is a
{\em block} of $M\in {\mathcal M}_{n,p}$ if, for every pair $(i,j)$
with $i\in A$ and $j\in \{1,\ldots,p\}\setminus A$, the matrix element
$M_{i,j}$ is minimal, and if no proper subset of $A$ has this property.
\end{defin}
\begin{rema}
The blocks of any matrix $M\in{\mathcal M}_{n,p}$ are disjoint, and there
is a unique decomposition of the full set of indices $\{1,\ldots,p\}$
into blocks.
\end{rema}
\begin{defin}A set of indices $A\subseteq\{1,\ldots,p\}$ is {\em contiguous}
if for any $i,k\in A$ and $j\in\{1,\ldots,n\}$, $i<j<k$ implies
$j\in A$.
\end{defin}
\begin{defin}
A matrix $M\in{\mathcal M}_{n,p}$ has {\em block form} if all of its blocks
are contiguous.
\end{defin}
We will also use the term ``block'' to refer to the sub-matrix
obtained by removing all rows and columns whose indices are not in the block.
\begin{defin}
The {\em partial row} $P_{r,k}(M)$ with $M\in {\mathcal M}_{n,p}$ and
$1\le r,k\le p$ is the $k$-dimensional vector formed by taking the first
$k$ components of row $r$ of $M$.
\end{defin}
\begin{notat} For integers $a$ and $b$ we use $a\la b$ to indicate that
$|a| < |b|$, with similar notation
for other inequality signs.  For integer vectors $\uvec$ and $\vvec$,
$\uvec\la\vvec$ indicates that $\uvec$ is lexicographically less than $\vvec$
where comparison of individual vector elements is made using $\la$.
\end{notat}
\begin{defin}
A matrix $M\in{\mathcal M}_{n,p}$ is {\em lexicographically ordered} if its
rows are in descending lexicographic order, according to the ordering $\la$.
\end{defin}
Since each row of $M$ is dominated by its diagonal element, this definition is
equivalent to the criterion that, for all pairs of indices $1\le i<j\le p$,
the partial rows satisfy $P_{i,i-1}(M) \ga P_{j,i-1}(M)$.
\begin{lem}
A lexicographically ordered matrix $M\in{\mathcal M}_{n,p}$, has block form.
\end{lem}
\begin{proof} Assume the contrary.  Then there is at least one pair of
blocks, $(A,B)$, whose indices are interleaved.  Assume the lowest index in
$A$ to be less
than the lowest index $j$ in $B$.  Then $A$ is partitioned into two subsets:
$A_{<j}$ of indices less than $j$ and $A_{>j}$ of indices greater than $j$.
There exists $i\in A_{<j}$ and $k\in A_{>j}$ such that $M_{i,k}$ is
non-minimal, since otherwise $A_{<j}$ and $A_{>j}$ would be separate blocks.

By the lexicographic ordering of $M$ we have $P_{j,j-1}(M)\gea P_{k,j-1}(M)$.
However, all entries of $P_{j,j-1}(M)$ are minimal since $j$ is the smallest
index in its block, whereas $M_{k,i}=M_{i,k}$ is non-minimal ---
a contradiction.
\end{proof}

\begin{notat}
Let $M_1\in{\mathcal M}_{n,p_1}$ and $M_2\in{\mathcal M}_{n,p_2}$
We say that $M_1\la M_2$ if either
\begin{enumerate}
\item there exists an index $k$ such that $1\le k\le\min(p_1,p_2)$,
$P_{i,i-1}(M_1)=P_{i,i-1}(M_2)$ for all $i<k$, and
$P_{k,k-1}(M_1)\la P_{k,k-1}(M_2)$, or
\item $p_1<p_2$ and $P_{i,i-1}(M_1)=P_{i,i-1}(M_2)$ for all $i\le p_1$.
\end{enumerate}
\end{notat}

Clearly the blocks of a lexicographically ordered matrix (considered as
sub-matrices) are themselves lexicographically ordered.

There may be many
lexicographically ordered matrices equivalent to a given matrix.
Because $\la$ is a total ordering and ${\mathcal M}_{n,p}$ is finite,
any subset of ${\mathcal M}_{n,p}$ has a greatest element.
\begin{defin}
A matrix $M\in{\mathcal M}_{n,p}$ is {\em lexicographically maximal} if it
is the greatest element of ${\mathcal E}_M$ under the ordering $\la$y.
\end{defin}
Trivially, a matrix that is lexicographically maximal is lexicographically
ordered.  The proof of the following is straightforward.
\begin{lem}
A matrix $M\in{\mathcal M}_{n,p}$ is lexicographically maximal if and only if
all its blocks (considered as sub-matrices) are lexicographically maximal
and are arranged in descending lexicographic order along the diagonal.
\end{lem}

In Section~\ref{sect:Gram}, we present a backtracking algorithm to generate
a list containing
a canonical representative of each equivalence class in
${\mathcal M}_n(\dmin)$.  The canonical representative of an equivalence
class is taken to be its lexicographically maximal element.  The algorithm
builds up matrices by augmenting the leading sub-matrix, one row and column
per iteration,
and is so designed that any sub-matrix $M\in{\mathcal M}_{n,p}$, constructed at
stage $p$, is automatically in block form with lexicographically ordered blocks
in descending lexicographic order along the diagonal.  To ensure that only
canonical matrices are produced in the end, there is a check, performed at the
completion of any block, that the block is lexicographically maximal.  In
the event that it is not, the sub-matrix is discarded.  To carry out this 
check, the following simple procedure, while not optimal, proved to be
adequate.

\subsection{Procedure {\tt IsLexMax}}\label{subsect:IsLexMax}
The procedure returns {\tt True} if the block cannot be transformed into some
lexicographically greater block by permutation of indices, and returns
{\tt False} otherwise.  It is perhaps simpler to consider how to
put a block into its lexicographically maximal order.  There is at least
one permutation that accomplishes this.  If the set
of all such permutations includes the identity,
then our routine is to return {\tt True}, otherwise {\tt False}.

Let the block be of size $p\ge2$.
The procedure goes as follows: we first find the set, ${\mathcal P}_2$, of
permutations that produce the lexicographically maximal $2\times2$ leading
sub-matrix.  Since this is
accomplished by placing a maximal off-diagonal element into the $(1,2)$
position, the set ${\mathcal P}_2$ can contain permutations of three types,
depending on the initial
locations of the maximal elements: the identity, if a maximal element
already was in the $(1,2)$ position; a single swap, if a maximal element
occurred in a position with an index (row or column) 1 or 2; a pair of swaps,
if a maximal element occurred in a position with neither index equal to 1 or
2.

Now, given the set ${\mathcal P}_k$ of permutations that produce the maximal
$k\times k$ leading sub-matrix, $2\le k\le p-1$, we produce the set
${\mathcal P}_{k+1}$.  We do this by composing each permutation in
${\mathcal P}_k$ with each single swap
$(k+1,j)$, $j>k+1$, that produces the maximal $(k+1)\times(k+1)$ leading
sub-matrix (if any exists), and with the identity, if the leading
$(k+1)\times(k+1)$
sub-matrix already is maximal.  The permutations so produced are the elements
of ${\mathcal P}_{k+1}$.  It is in fact unnecessary to consider all
possible swaps, since some indices in the set $\{k+2,\ldots,p\}$ may be
equivalent (under permutation of indices) to each other due to the
structure of the matrix.  Thus
we need only consider one index from each equivalence class.

This process is iterated until either
\begin{enumerate}
\item a set ${\mathcal P}_k$ does not contain the
identity, in which case we terminate the procedure
and return {\tt False}, or
\item we have produced ${\mathcal P}_p$.  If it contains
the identity, we return {\tt True}, otherwise {\tt False}.
\end{enumerate}
Notice that multiplying all elements of the set ${\mathcal P}_p$ by the inverse
of any one of them produces the automorphism group of the block.

\section{Finding candidate Gram matrices}\label{sect:Gram}
We perform a backtracking search, based on the methods of Moyssiadis and
Kounias~\cite{MK} and Chadjipantelis, Kounias and Moyssiadis~\cite{CKM}, to
determine the set ${\mathcal M}_n^\square(\dmin)$.  Starting from the
$1\times1$ matrix, $\begin{pmatrix}n\end{pmatrix}$, candidate matrices are
built up by symmetrically appending one row and column at a time,
until size $n$ is reached or no continuation is possible.
At this stage, the program
returns to the most recent sub-matrix from which there is a possible
continuation that has not yet been tried, and resumes the search from there.
Proceeding exhaustively in this
manner, until the entire search tree has been explored, the program
terminates with a complete list of candidate Gram matrices.

Naturally, the search tree is vast, and various methods must
be used to prune it.  
The value $\dmin^2$  bounds the space we need to explore.  In practice,
a theorem of Moyssiadis and Kounias, described below, provides a criterion
for removing branches of the search tree at an early stage.
However, the size of the search space grows rapidly as $\dmin$ is lowered,
so there is a practical limit to how low the threshold can be set.  

A second method to reduce the size of the search is to impose
lexicographic maximality of the matrices.  This means that all but one
of the branches containing a set of equivalent matrices are removed.
There are two aspects to this pruning.  The first
is that only lexicographically ordered matrices whose blocks are
in descending lexicographic order along the diagonal are generated.
Specifically, there is an {\em active block}, namely the block
containing the index of the most recently added row and column.  A new
block begins in index $j$ whenever $P_{j,j-1}(M)$ is a minimal vector.
The row and column vector appended to the current sub-matrix must have
non-minimal elements only within the active block, and must be lexicographically
less than the preceding row and column.  (The first criterion is actually
redundant.) Furthermore, a comparison is made of the active block 
(considered as a sub-matrix) and the previous block, and the continuation
is disallowed if it would cause the blocks to be incorrectly ordered.

However, as discussed in
Section~\ref{subsect:GramProperties}, this alone does
not guarantee that the matrix obtained is lexicographically
maximal.  To enforce maximality, the active
block is subjected to the test {\tt IsLexMax}, described in
Section~\ref{subsect:IsLexMax}, and the continuation is disallowed if
the test fails.
This test can be extremely time-consuming, although it improved the efficiency
of the algorithm for the cases considered here.  For higher orders, it may
be better to omit it, or, better, to carry it out more efficiently, perhaps
using the methods of~\cite{McK}.

\subsection{The theorem of Moyssiadis and Kounias}\label{subsect:MKTheorem}
\begin{theo}\label{MKTheorem}
Let $M=\begin{pmatrix}D_r & B\\B^{\text{T}} & A\end{pmatrix}$ be a 
symmetric, positive definite matrix of order $m$ with elements taken
from the set $\Phi$ whose members are greater or equal in magnitude to some
positive number $c$.  Here $D_r$ and $A$ are square matrices and the
order of $D_r$ is $r\le m\le n$.  The diagonal elements of $A$ are equal
to $n$.  The columns of the $r\times(m-r)$ matrix $B$ are taken from some set
$\Gamma_r\subseteq \Phi^r$.

Define $d^*$ and $\gamma^*$ by
\begin{equation}
d^*=\begin{vmatrix}D_r & \gamma^*\\ \gamma^{*\text{T}} & c \end{vmatrix}
=\underset{\gamma\in\Gamma_r}{\max}
\begin{vmatrix}D_r & \gamma\\ \gamma^{\text{T}} & c \end{vmatrix}
\end{equation}
and define
\begin{equation}\label{MKboundExpr}
u_r(d):=(n-c)^{m-r-1}\left[(n-c)\det D_r+(m-r)\max(0,d)\right].
\end{equation}
Then
\begin{equation}
\label{MKbound}
\det M\le u_r(d^*)
\end{equation}
\end{theo}
For the convenience of the reader, we will reproduce the proof of Moyssiadis
and Kounias below.  First we point out two slight modifications of their
original statement of the theorem.
\begin{itemize}
\item The columns of $B$ and the
vector $\gamma$ are taken to be members of $\Gamma_r\subseteq\Phi^r$
rather than of $\Phi^r$ itself.  This is needed in a practical
implementation of the algorithm where lexicographic ordering is enforced.
\item We do not require $d^*\ge0$ but replace $d$ by
$\max(0,d)$ in the bound~(\ref{MKboundExpr}).  This is necessary in the
case $n\equiv3\bmod4$ since, in contrast to the case $n\equiv1\bmod4$,
$d^*$ is often negative.
\end{itemize}
\begin{proof}
By subdivision of the matrices $A$ and $B$ we rewrite $M$ as
\begin{equation*}
\begin{pmatrix}
D_r & B_1 & \gamma\\
B_1^{\rm T} & A_1 & \delta\\
\gamma^{\rm T} & \delta^{\rm T} & n
\end{pmatrix}
\end{equation*}
where $\gamma\in\Gamma_r$ and $\delta$ is a column vector of length $m-r-1$.
Using the expansion by minors on the last column of $M$ we obtain
\begin{equation*}\label{MKinduct}
\det M=(n-c)\begin{vmatrix}
D_r & B_1\\
B_1^{\rm T} & A_1
\end{vmatrix}
+\begin{vmatrix}
D_r & B_1 & \gamma\\
B_1^{\rm T} & A_1 & \delta\\
\gamma^{\rm T} & \delta^{\rm T} & c
\end{vmatrix}.
\end{equation*}
Note that the first term is an $(m-1)\times(m-1)$ determinant that
satisfies the hypotheses of the theorem.  Denote the array in the second
term by $\Mbar$.  Either $\det\Mbar\le0$ or $\Mbar$ is positive definite
in which case we have
\begin{equation*}
\det\Mbar=c\begin{vmatrix}
D_r-\frac{1}{c}\gamma\gamma^{\rm T} & B_1-\frac{1}{c}\gamma\delta^{\rm T} \\
B_1^{\rm T}-\frac{1}{c}\delta\gamma^{\rm T} & A_1-\frac{1}{c}\delta\delta^{\rm T}
\end{vmatrix}
\le c\begin{vmatrix}
D_r-\frac{1}{c}\gamma\gamma^{\rm T}
\end{vmatrix}
\cdot\begin{vmatrix}
A_1-\frac{1}{c}\delta\delta^{\rm T}
\end{vmatrix}.
\end{equation*}
The first factor is bounded above by $d^*$, and the second,
being the determinant of a positive definite matrix whose diagonal entries
are less than or equal to $n-c$, is bounded above by $(n-c)^{m-r-1}$.
Hence
\begin{equation*}\label{MbarBound}
\det\Mbar\le(n-c)^{m-r-1}\max(0,d^*).
\end{equation*}
By induction, using this bound and~(\ref{MKinduct}) we obtain~(\ref{MKbound}).
\end{proof}

\begin{coro}\label{rZeroCoro}
Let $M$ be an $m\times m$ symmetric, positive definite matrix whose diagonal
elements are less than or equal to some positive number $n$ and whose
off-diagonal
elements are greater than or equal to some positive number $c$ in magnitude.
Then
\begin{equation*}
\det M\le (n-c)^m+mc(n-c)^{m-1}
\end{equation*}
\end{coro}
\begin{proof}
Since $M$ is positive definite, its determinant can only increase if
we replace all diagonal elements with the upper bound, $n$.
Now apply Theorem~(\ref{MKTheorem}) with $r=0$ to this revised matrix.
\end{proof}

In the $n\equiv1\bmod4$ case our search algorithm will make use of
Theorem~(\ref{MKTheorem}) directly (with $c=1$ and $m=n$).
In the $n\equiv3\bmod4$ case we strengthen the theorem as follows:
\begin{coro}\label{nequals3Coro}
Let $M=\begin{pmatrix}D_r & B\\B^{\text{T}} & A\end{pmatrix}$ satisfy
the conditions of Theorem~(\ref{MKTheorem}) with $m=n\equiv3\bmod4$ and
furthermore let $M\in\mathcal M_{n,n}$.  Then
\begin{equation*}
\det M\le(n-1)^{n-r}\det D_r + \left[(n-1)^{n-r}-(n-3)^{n-r}
-(n-r)(n-3)^{n-r-1}\right]\max(0,d^*).
\end{equation*}
\end{coro}
\begin{proof}
We apply Corollary~\ref{rZeroCoro} to the matrix
$A'_1:=A_1-\frac{1}{c}\delta\delta^{\rm T}$ in the proof of
Theorem~\ref{MKTheorem}
to obtain a bound on $\det\Mbar$ which is sharper than~(\ref{MbarBound}).
Because
all elements of both $A_1$ and $\delta$ are congruent to 3 modulo 4, it
follows that the diagonal elements of $A'_1$ are bounded above by $n-1$
and the magnitude of the off-diagonal elements is bounded below by $2$.
Therefore
\begin{equation*}
\det A'_1\le(n-3)^{m-r-1}+2(m-r-1)(n-3)^{m-r-2}
\end{equation*}
and $\det\Mbar$ is bounded above by the product of this bound and
$\max(0,d^*)$.  Applying the induction step, we find that
the coefficient of $\max(0,d^*)$ in the bound on $\det M$ becomes the sum
\begin{equation*}
\sum_{j=0}^{n-r-1}(n - 1)^{n-r-1-j} [(n - 3)^j + 2j(n - 3)^{j-1}].
\end{equation*}
Evaluating this sum gives the stated result.
\end{proof}

\subsection{The algorithm}
Let the set of allowed matrix elements be
\begin{equation*}
\Phi=\begin{cases}
\{1,-3,\ldots,2-n\} & \text{if $n\equiv1\bmod4$} \\
\{-1,3,\ldots,2-n\} & \text {if $n\equiv3\bmod4$}
\end{cases}
\end{equation*}
and denote the determinant threshold (for $\{-1,1\}$-matrices) as $\dmin$.
The threshold for candidate Gram matrices is then $\dmin^2$.

We carry out the following recursive procedure.
\begin{enumerate}
\item Initialize variables.
\begin{align*}
r & := 1 && \text{order of current sub-matrix} \\
M_1 & := \begin{pmatrix}n\end{pmatrix} && \text{the initial sub-matrix} \\
F_0 &:= (()) && \text{the list of vectors (Here $()$ is the null vector.)} \\
\ebar_1 &:= n-2 && \text{the magnitude of the maximum allowed matrix element}
\end{align*}
\item (Beginning of the recursive step)
Increment $r$.  Build the list of vectors, $\Ftilde_{r-1}$, by
appending to every vector in the list
$F_{r-2}$ each $e\in\Phi$ satisfying $|e|\le\emax$.
The construction is carried out in such a way that $\Ftilde_{r-1}$ will
always be in ascending lexicographic order (according to the ordering $\la$).
Initialize $F_{r-1}$ to the null list.  Build
the list $\Gamma_r$ by appending to every vector in $F_{r-2}$ all possible
2-vectors, $(e_1,e_2)\in\Phi^2$ satisfying $|e_1|,|e_2|\le\ebar_{r-1}$.
\item For each vector $f$ in $\Ftilde_{r-1}$ construct the matrix
\begin{equation*}
M_r=\begin{pmatrix}
M_{r-1} & f\\
f\T & n
\end{pmatrix}.
\end{equation*}

If $r=n$ then
\begin{enumerate}
\item If $\det M_r\ge\dmin^2$ and is a perfect square, then run
the procedure {\tt IsLexMax} on the final block of $M_r$.  If it returns
{\tt True} then print out the matrix.
\item Return to the calling procedure.
\end{enumerate}
If $r<n$ go to step (4).
\item For each $\gamma\in\Gamma_r$ evaluate
\begin{equation*}
d=\begin{vmatrix}
M_r & \gamma \\
\gamma\T & 1
\end{vmatrix}
\end{equation*}
until a $d$ is found such that the bound in Theorem~\ref{MKTheorem} or
Corollary~\ref{nequals3Coro} equals or exceeds $\dmin^2$, or until the list
$\Gamma_r$ is exhausted.  If
the list was exhausted without finding such a $d$, then
$M_r$ is discarded, and we continue in step (3) with the next vector $f$.
Otherwise, add $f$ to the list $F_{r-1}$.  Then test the active block
for lexicographic maximality using {\tt IsLexMax}.  If it passes,
recursively begin again at step (2), otherwise continue in step (3) with
the next vector $f$.
\end{enumerate}

{\bf Notes:}
\renewcommand{\labelenumi}{(\alph{enumi})}
\begin{enumerate}
\item If a new block started with index $r$, then index $r+1$ is
very important.  It must be the case that all components of the
vector $f\in\Ftilde_r$ that forms the initial segment of row and
column $r+1$ are minimal except possibly the
$r^{\text{th}}$.  If this component is minimal as well, then the entire
remainder of the matrix $M_n$ will be filled in with minimal elements, and
we can proceed directly to the check of the perfect-square condition as in
step (3).
If component $r$ is non-minimal, we set $\ebar_r$ equal to its absolute value.
\item No attempt has been made to improve the efficiency of {\tt IsLexMax}
by keeping track of the automorphisms of the blocks of $D_r$ as they are
built.  This may be necessary to extend to method to higher orders in the
$n\equiv3\bmod4$ case.
\item No attempt has been made to improve the efficiency of step (4) by
using facts about the form of the optimal form of $\gamma$.  Moyssiadis
and Kounias~\cite{MK}, for example, in their analysis of $n=17$ were able
to predict
when a component would equal $1$ which avoided the need to check the entire
list $\Gamma_r$.  This condition occurs much less frequently in the case
$n\equiv3\bmod4$.
\end{enumerate}
\renewcommand{\labelenumi}{(\arabic{enumi})}

\section{Order by order list of candidate Gram matrices}\label{sect:Mresults}
We state known results for all odd orders up to $n=21$ and for selected
higher orders.  First we introduce
a few standard types of matrices that arise repeatedly.
\begin{notat} The matrix $J_{p,q}$ is the $p\times q$ matrix with all elements
equal to 1.  We define $J_p:=J_{p,p}$.
\end{notat}
Barba's bound is attained by the determinant of a $\{-1,1\}$-matrix, $R$, if
and only if $M_r(R)=M_c(R)=(n-1)I_n+J_n$.  This occurs in orders 1, 5, 13,
25, and 41, among others, and we will have nothing more to say about these
orders.  No other candidate Gram matrix can have a determinant as large
as Barba's bound, but when $n\equiv1\bmod4$ it is not unusual to find candidate
Gram matrices of high determinant which differ only slightly from the the
matrix $(n-1)I_n+J_n$.  Typically the off-diagonal entries differ from 1 in
only one or two rows and columns.
\begin{notat} If $v$ is a vector of length $n-1$ then $S(v)$ denotes
the matrix
\begin{equation*}
\begin{pmatrix}
n &  v^{\rm T}\\
v & (n-1)I_{n-1}+J_{n-1}
\end{pmatrix}
\end{equation*}
\end{notat}
\begin{notat} If $v$ and $w$ are vectors of length $n-2$ then $D(a;v;w)$
denotes the matrix
\begin{equation*}
\begin{pmatrix}
n & a & v^{\rm T} \\
a & n & w^{\rm T} \\
v & w & (n-1)I_{n-2}+J_{n-2}
\end{pmatrix}
\end{equation*}
\end{notat}

For $n\equiv3\bmod4$ define
\begin{defin}
An {\em Ehlich block matrix} is a matrix of the form
\begin{equation*}
B(v_1,v_2,\ldots,v_k):=(n-3)I_{n}-J_n+4
\begin{pmatrix}
J_{v_1} & 0 & \ldots & 0  \\
0 & J_{v_2} & & \\
\vdots & &\ddots & & \\
0 & 0 & \ldots & J_{v_k} \\
\end{pmatrix}
\end{equation*}
with $\sum v_j=n$, or in other words, a matrix with diagonal elements equal
to $n$, square blocks along the diagonal whose off-diagonal elements are
equal to 3, and all other off-diagonal elements $-1$.
\end{defin}
In low orders, $n\equiv3\bmod4$, most large-determinant candidate Gram matrices
are closely related to Ehlich block matrices.
When an element $a$ of a vector is repeated $k$ times consecutively, it will be
convenient to denote this sequence by $a_k$.  

We describe below the candidate Gram matrices found by our procedure or, in
a few cases which we note explicitly, matrices which were found by other
means.  In higher orders, where the the number of matrices tends to be
large, we sometimes omit the explicit listing of matrices.  Complete lists
of matrices are available from the author upon request.

\subsection{$n=3$}
The largest determinant is $1\times2^2$.  Our program finds the unique
candidate Gram matrix, $B(1_3)$, corresponding to this value, and no
other candidates of equal or larger determinant.

\subsection{$n=7$}
The maximal determinant, found by Williamson~\cite{Wi}, is $9\times2^6$.
It corresponds to the matrix $B(2_3,1)$, which is found by our program.
No other candidates of equal or larger determinant are found.

\subsection{$n=9$}
The maximal determinant, found by Ehlich and Zeller~\cite{EZ}, is
$56\times2^8$.  It corresponds to the candidate Gram matrix $S(5,1_7)$.
No other candidates with equal or larger determinant are found.

\subsection{$n=11$}
This case was treated by Ehlich, who found that the maximal determinant
is $320\times2^{10}$, but his analysis was never published.  His
results were reported in the paper of Galil and Kiefer~\cite{GK}.  There
are seven candidate
Gram matrices.  Four of these have determinant $(324\times2^{10})^2$, above
the value corresponding to the maximal determinant.
They are $B(3,2,1_6)$, $B(4,1_7)$, $B(5,1_6)$ and
\begin{equation*}
\begin{pmatrix}
11&3&3&&&\\
3&11&-1&&-J_{3,8}&\\
3&-1&11&&&\\
&&&&&\\
&\mspace{-18.0mu}-J_{8,3}\mspace{-18.0mu}&&&B(4,1_4)&\\
&&&&&&&\\
\end{pmatrix}.
\end{equation*}
The failure of these four matrices to produce a maximal determinant is
discussed in the next section.

The other three candidates all correspond to the maximal determinant.
They are $B(5,2_3)$,
\begin{equation*}
\begin{pmatrix}
11&3&3&-&-&&&\\
3&11&-&3&-&&&\\
3&-&11&-&3&&-J_{5,6}&\\
-&3&-&11&-&&&\\
-&-&3&-&11&&&\\
&&&&&&&\\
&&\mspace{-18.0mu}-J_{6,5}\mspace{-18.0mu}&&&&B(2_3)&\\
&&&&&&&\\
\end{pmatrix}, \text{ and }
\begin{pmatrix}
11&3&3&3&-&&&\\
3&11&3&-&-&&&\\
3&3&11&-&-&&-J_{5,6}&\\
3&-&-&11&3&&&\\
-&-&-&3&11&&&\\
&&&&&&&\\
&&\mspace{-18.0mu}-J_{6,5}\mspace{-18.0mu}&&&&B(2_3)&\\
&&&&&&&\\
\end{pmatrix},
\end{equation*}
where ``$-$'' stands for $-1$.  All candidates have unique characteristic
equations.

\subsection{$n=15$}
All four candidate Gram matrices are Ehlich block matrices.
The matrices $B(4_3,3)$ and $B(6,3,2_3)$ have determinant
$(105\cdot3^5\cdot2^{14})^2$, corresponding to the conjectured maximal
value, while $B(3_4,2,1)$ and $B(3_5)$ have
determinant $(108\cdot3^5\cdot2^{14})^2$.  All of the matrices have
different characteristic equations.
Our algorithm, implemented as a Mathematica program, was able
to generate the four candidates in about seven hours on a personal computer.
This is a far larger
running time than for any other order for which $\md(n)$ had
been established, all of which
require anywhere from a few milliseconds to a few minutes to produce the
complete list of candidates.

\subsection{$n=17$}
The maximal determinant for $n=17$ was conjectured by Schmidt~\cite{sch}
to be $5\times4^8\times2^{16}$.  The conjecture was independently
formulated and proved by Moyssiadis and Kounias~\cite{MK} 
The candidate Gram matrices are $S(-3,-3,1_{14})$ with determinant
$(22\times4^7\times2^{16})^2$, $S(5,-3_2,1_{13})$ with determinant
$(21\times4^7\times2^{16})^2$, $D(-3;-3_3,1_{12};1_3,-3,1_{11})$ with
determinant $(83\times4^6\times2^{16})^2$, $D(1;5,-3,1_{13};1,-3_3,1_{11})$
with determinant $(81\times4^6\times2^{16})^2$, the matrix
\begin{equation*}
\begin{pmatrix}
  A&J_2&J_2&J_2&J_2& &         &\\
J_2&  B&J_2&J_2&J_2& &         &\\
J_2&J_2&  B&J_2&J_2& &         &\\
J_2&J_2&J_2&  B&J_2& &J_{10,7}&\\
J_2&J_2&J_2&J_2&  B& &         &\\
&&&&&&&&&\\
&&\mspace{-18.0mu}J_{7,10}\mspace{-18.0mu}&&&&16I_7+J_7&\\
&&&&&&&&&\\
\end{pmatrix}\text{ where }
A=\begin{pmatrix} 17 & 5\\ 5 & 17\end{pmatrix},\quad
B=\begin{pmatrix} 17 & -3\\ -3 & 17\end{pmatrix}
\end{equation*}
with determinant $(5175\times4^3\times2^{16})^2$, and the six matrices
\begin{align*}
M^{(17)}_1&=D(1;1,-3_6,1_8;-3,1_{14}),\\
M^{(17)}_2&=S(-3_8,1_8),\\
M^{(17)}_3&=S(-3_{16}),\\
M^{(17)}_4&=D(-3;-3_2,1_{13};-3_2,1_{13}),\\
M^{(17)}_5&=S(5_2,-3_2,1_{12}),\\
M^{(17)}_6&=
\begin{pmatrix}
17&-3&-3&-3& 1& 1&&         &\\
-3&17&-3& 1& 1& 1&&         &\\
-3&-3&17& 1& 1& 1&&         &\\
-3& 1& 1&17&-3&-3&&J_{6,11}&\\
 1& 1& 1&-3&17& 1&&         &\\
 1& 1& 1&-3& 1&17&&         &\\
&&&&&&&&&\\
&&&\mspace{-18.0mu}J_{11,6}\mspace{-18.0mu}&&&&16I_{11}+J_{11}&\\
&&&&&&&&&\\
\end{pmatrix}
\end{align*}
each with determinant $(5\times4^8\times2^{16})^2$.

The maximal determinant corresponds to the matrix $M^{(17)}_3$
For some reason, the two matrices $M^{(17)}_1$ and $M^{(17)}_5$
were omitted from the list in~\cite{MK}.  The characteristic
equations of the above
matrices are distinct, except for those of $M^{(17)}_2$ and $M^{(17)}_5$.

\subsection{$n=19$}
Our program fails to run to completion in a reasonable amount of time, or
to turn up any candidate matrices,
for any threshold below Ehlich's bound~(\ref{Bthree}).  Two Gram matrices,
corresponding to the best known determinant value, $833\times4^6\times2^{18}$,
are found by direct computation from the conjectured maximal $\{-1,1\}$
matrices of Smith~\cite{Sm} and Cohn~\cite{Co2}. 
They are
\begin{equation*}
\begin{pmatrix}
19& 3& 3& 3& -& -& -&&         &\\
 3&19& 3& 3& -& -& -&&         &\\
 3& 3&19& 3& -& -& -&&         &\\
 3& 3& 3&19& 3& 3& 3&&-J_{7,12}&\\
 -& -& -& 3&19& 3& 3&&         &\\
 -& -& -& 3& 3&19& 3&&         &\\
 -& -& -& 3& 3& 3&19&&         &\\
&&&&&&&&&\\
&&&\mspace{-18.0mu}-J_{12,7}\mspace{-18.0mu}&&&&&B(3_4)&\\
&&&&&&&&&\\
\end{pmatrix}
\end{equation*}
and
\setcounter{MaxMatrixCols}{15}
\begin{equation*}
\begin{pmatrix}
19& -& -& 3& -& -& 3& -& -& 3&&         &\\
 -&19& 3& 3& -& -& -& -& -& -&&         &\\
 -& 3&19& 3& -& -& -& -& -& -&&         &\\
 3& 3& 3&19& -& -& -& -& -& -&&         &\\
 -& -& -& -&19& 3& 3& -& -& -&&-J_{10,9}&\\
 -& -& -& -& 3&19& 3& -& -& -&&         &\\
 3& -& -& -& 3& 3&19& -& -& -&&         &\\
 -& -& -& -& -& -& -&19& 3& 3&&         &\\
 -& -& -& -& -& -& -& 3&19& 3&&         &\\
 3& -& -& -& -& -& -& 3& 3&19&&         &\\
&&&&&&&&&&&&\\
&&&&\mspace{-18.0mu}-J_{9,10}\mspace{-18.0mu}&&&&&&&B(3_3)&\\
&&&&&&&&&&&&\\
\end{pmatrix}.
\end{equation*}
\setcounter{MaxMatrixCols}{10}
Extensive searching by various random methods has failed to improve upon
this determinant value.
A simple exhaustive search for Ehlich block matrices with larger perfect square
determinant than that of the above matrices yields the result that no
such matrices exist.

\subsection{$n=21$}
The largest determinant for $n=21$ is $29\times5^9\times2^{20}$ which was found
by Chadjipantelis, Kounias, and Moyssiadis~\cite{CKM} and proved by
them to be maximal.
The candidate Gram matrices are $S(-3_5,1_{15})$ with determinant
$(6\times5^{10}\times2^{20})^2$, the matrix
\begin{equation*}
\begin{pmatrix}
  A&J_2&J_2&J_2&J_2& &         &\\
J_2&  B&J_2&J_2&J_2& &         &\\
J_2&J_2&  B&J_2&J_2& &         &\\
J_2&J_2&J_2&  B&J_2& &J_{10,11}&\\
J_2&J_2&J_2&J_2&  B& &         &\\
&&&&&&&&&\\
&&\mspace{-18.0mu}J_{11,10}\mspace{-18.0mu}&&&&16I_{11}+J_{11}&\\
&&&&&&&&&\\
\end{pmatrix}\text{ where }
A=\begin{pmatrix} 21 & 5\\ 5 & 21\end{pmatrix},\quad
B=\begin{pmatrix} 21 & -3\\ -3 & 21\end{pmatrix}
\end{equation*}
with determinant $(18432\times5^5\times2^{20})^2$, and the three matrices
$M^{(21)}_1=S(5,-3_5,1_{14})$, $M^{(21)}_2=S(5^4,1_{16})$, and
$M^{(21)}_3=S(-7,-3_2,1_{17})$ all with determinant
$(29\times5^9\times2^{20})^2$.  $M^{(21)}_1$ and $M^{(21)}_3$ have the
same characteristic equation; all others are distinct.  The matrix
corresponding to the maximal determinant is $M^{(21)}_2$.  For some reason
the matrices $M^{(21)}_1$ and $M^{(21)}_3$ were omitted from the list
in~\cite{CKM}.

\subsection{$n=29$, $33$, and $37$}
In orders $29$ and $33$ the largest known determinants were found by
Bruce Solomon using the gradient ascent algorithm described in~\cite{OSDS}.
In order $37$ the largest known determinant was constructed in~\cite{OS}.

Unfortunately, in all three cases these lower bounds are too low to use as
thresholds in the search for candidate Gram matrices.  Both time and space
requirements appear to be growing exponentially as the threshold is lowered.
The best we have been able to do is narrow the range within which $\md(n)$
must lie, by using somewhat higher thresholds.  The thresholds used, and
the number of candidates found for each determinant value equal to or
exceeding the given threshold, are shown in Table~\ref{table:29-33-37}.
\begin{table}
\begin{tabular}{cc|cc|cc}
$n=29$ & & $n=33$ & & $n=37$\\
\hline
det. & num. & det. & num. & det. & num.\\
\hline
$51\times7^{13}$ & 2 & $495\times8^{14}$ & 1 & $680\times9^{16}$ & 1\\
$355\times7^{12}$ & 1 & $252315\times8^{11}$ & 1 & $75\times9^{17}$ & 1,1\\
$352\times7^{12}$ & 1 & $3929\times8^{13}$ & 1 & $672\times9^{16}$ & 1\\
$50\times7^{13}$ & 1,2 & $490\times8^{14}$ & 1,1 & $485070\times9^{13}$ & 1\\
$2448\times7^{11}$ & 1 & $3919\times8^{13}$ & 1 & $665\times9^{16}$ & 1\\
$119808\times7^9$ & 1 & $489\times8^{14}$ & 1 & $53760\times9^{14}$ & 1\\
$349\times7^{12}$ & 1 & $3911\times8^{13}$ & 1 & $4352000\times9^{12}$ & 1\\
$348\times7^{12}$ & 1 & $250047\times8^{11}$ & 1 & $663\times9^{16}$ & 1,1\\
$2432\times7^{11}$ & 1 & $1023942465\times8^7$ & 1 & $483000\times9^{13}$ & 1\\
$2430\times7^{11}$ & 1 & $3906\times8^{13}$ & 3 & $4345600\times9^{12}$ & 1\\
$347\times7^{12}$ & 2 & $61\times8^{15}$ & $1_3$ & $5952\times9^{15}$ & 1\\
$118656\times7^9$ & 1 & $31203\times8^{12}$ & 1 & $661\times9^{16}$ & 2\\
$2416\times7^{11}$ & 1 & $3898\times8^{13}$ & 2 & $5946\times9^{15}$ & 1\\
$345\times7^{12}$ & 1,3,4 & $3897\times8^{13}$ & 1,1 & $53504\times9^{14}$ & 1\\
$2413\times7^{11}$ & 1 & $3895\times8^{13}$ & 1 & $5942\times9^{15}$ & 2\\
$118128\times7^9$ & 1 & $31131\times8^{12}$ & 1 & $660\times9^{16}$ & $1_3$\\
$344\times7^{12}$ & $1_3$ & $3889\times8^{13}$ & 1,2 & $659\times9^{16}$ & 1\\
$2403\times7^{11}$ & 1 & $31108\times8^{12}$ & 1\\
$49\times7^{13}$ & $1_4$,2 & $486\times8^{14}$ & $1_5$,3\\
$2400\times7^{11}$ & $1_6$,2 & $31098\times8^{12}$ & 1\\
$117504\times7^9$ & $1_3$ & $3887\times8^{13}$ & 1\\
$2397\times7^{11}$ & 3 & $248724\times8^{11}$ & 1\\
$5750784\times7^7$ & 1 & $3885\times8^{13}$ & 1,2\\
$342\times7^{12}$ & $1_3$,2,3,7 & $31050\times8^{12}$ & 4\\
&& $3881\times8^{13}$ & 1,1\\
\end{tabular}
\caption{Exhaustive enumeration of candidate Gram matrices of orders 29,
33, and 37 with thresholds $342\times7^{12}$, $485\times8^{14}$, and
$659\times9^{16}$.  Square roots of the determinants of the candidates,
with the factor $2^{n-1}$ omitted, are given.  In
number of matrices column, multiple entries refer to multiple
characteristic equation classes.\label{table:29-33-37}}
\end{table}

\section{Decomposability of candidate Gram matrices}\label{sect:decomposition}

\subsection{General considerations}

Given a candidate Gram matrix, $M$, we would like to determine whether 
it admits a decomposition $M=RR^{\rm T}$ where $R$ is a $\{-1,1\}$-matrix.
This task appears to be rather difficult when only the matrix $M$ is given.
In the present context however, more information is available, namely the
complete list of candidate Gram matrices with determinant
equal to a specified perfect square value.  In many cases one can
demonstrate the nonexistence of decompositions using this knowledge.

Let $M_r$ be a candidate Gram matrix from our list, and suppose that
$M_r=RR^{\rm T}$ for some $\{-1,1\}$-matrix $R$.   Then, by
Lemma~\ref{lemma:charEqn} there is a matrix, $M_c=R^{\rm T}R$,
which has the same determinant and characteristic polynomial as $M_r$  Hence
by permuting columns of $R$ we can convert $M_c$ to an equivalent matrix 
which also appears on
our list.  A first step in our procedure is therefore to enumerate the
pairs $(M_r,M_c)$ of
candidate Gram matrices taken from our last which have the same characteristic
polynomial.  (Naturally, for every $M$ on our list, $(M,M)$ will be such a
pair, but in general there may be other pairs too.)

The pair of formulas
\begin{equation}\label{decompPair}
M_r=RR^{\rm T}, \qquad M_c=R^{\rm T}R
\end{equation}
implies a number of constraints on the possible rows and columns of $R$,
\begin{align}
& M_r^2 = RM_cR^{\rm T}\label{quad}\\
& M_c^2 = R^{\rm T}M_rR\label{quadT}\\
& M_rR = RM_c\label{rcconsistmat},
\end{align}
which are verified by substituting for $M_r$ and $M_c$
using~(\ref{decompPair}).
Let $x_j^{\rm T}$ denote the $j^{\text{th}}$ row of $R$, expressed as
a row matrix, and let $y_j$ denote the $j^{\text{th}}$ column of $R$,
expressed as a column matrix.  The rows of $M_r$ are
denoted by $m_j^{\rm T}$, while those of $M_c$ are denoted by ${m'}_j^{\rm T}$.
Using the notation $x_{ji}$ for the $i^{\text{th}}$ element of the
column matrix $x_j$, we have $x_{ji}=y_{ij}$.  The diagonal elements
of~(\ref{quad}) give a constraint on the rows of $R$,
\begin{equation}\label{rowquad}
m_j^{\rm T}m_j=x_j^{\rm T}M_c x_j
\end{equation}
while those of~(\ref{quadT}) give a constraint on the columns,
\begin{equation}\label{colquad}
{m'}_j^{\rm T}m'_j=y_j^{\rm T}M_r y_j.
\end{equation}

Finding solutions to these equations is facilitated by the large degree
of symmetry that most candidate Gram matrices with large determinant have.
A permutation of
rows, represented by a permutation matrix $P$, is an {\em automorphism} of
$M$ if $PMP^{\rm T}=M$.  (Multiplying on the right by $P^{\rm T}$ applies
the same permutation to columns, preserving the symmetry of $M$.)  The
set of automorphisms forms a group under composition, the {\em automorphism
group} of $M$.

A special type of permutation, the {\em swap}, interchanges a pair of rows.
The subset of the automorphism group generated by automorphisms which
are swaps forms a subgroup of the full automorphism group.  The
permutation that interchanges rows $i$ and $j$
is written $(ij)$.  If $(ij)$ and $(ik)$ are automorphisms of $M$, then
$(jk)$ is as well.  Hence we may order the rows of $M$ so that the rows
$j$ such that $(ij)$ is an automorphism, plus row $i$, form a contiguous block.
The lexicographic ordering we imposed earlier already implies
such a structure.  It is then easy to see that the automorphism subgroup
consisting of swaps is a direct sum of symmetric groups acting on these blocks.

When such structure is present, the right-hand side of~(\ref{rowquad})
will depend only on certain sums of elements of $x_j$, rather than on
the elements individually.  This greatly reduces the amount of work that
needs to be done, as should become clear in the examples that
follow.  Note that we only generate solutions which are
compatible with
parity normalization (see Section~\ref{sect:matProps}), that is, we require
$n+\sum_i x_{ji}\equiv 0 \bmod4$.  Similar considerations apply to the equation
for columns~(\ref{colquad}).


Once the row and column solutions have been found, we use the
off-diagonal elements of~(\ref{quad}) to determine which combinations
of row (or column) solutions are consistent.  The typical off-diagonal
constraint takes the form
\begin{equation}\label{rowconsist}
m_j^{\rm T}m_k=x_j^{\rm T}M_c x_k,
\end{equation}
and involves the sums of elements of $x_j$ mentioned above, as above, as well
as the inner products of the sub-vectors of rows $j$ and $k$ corresponding
to these sums.  The analysis is most straightforward when
the weights of these inner products are all equal.  Then only the inner
product of the full row, $x_j^{\rm T}x_k=m_{jk}$, enters rather than
the individual inner products.  We will find that this
occurs in all the examples studied below.


There is also a compatibility condition between row and column solutions,
derived from~(\ref{rcconsistmat}), which is sometimes useful, namely
\begin{equation}\label{rcconsist}
m_j^{\rm T}y_k=x_j^{\rm T}m'_k.
\end{equation}
In general the constraint will depend on the element $y_{kj}=x_{jk}$ and
on the sums of elements referred to previously.  In cases where the
off-diagonal elements of the diagonal block of $M_r$ containing $m_{jj}$
and of the diagonal block of $M_c$ containing $m'_{kk}$ are equal, the
dependence on $y_{kj}=x_{jk}$ is eliminated.  In the examples explicitly
worked out below, this constraint was not needed.

The above constraints, sometimes combined with other elementary consequences
of equation~(\ref{decompPair}), can often be used to prove the
non-decomposability
of a given pair of candidate Gram matrices.  We will illustrate the method
on the candidate Gram matrices for order 15, and also on a certain pair of
candidate Gram matrices of order 17.

\subsection{Order 15}

The four candidate Gram matrices listed in Section~\ref{sect:Mresults}
correspond to different characteristic polynomials.
Thus we must investigate four pairs
\begin{align*}
&(B(4_3,3),B(4_3,3)), &(B(6,3,2_3),B(6,3,2_3)),\\
&(B(3_5),B(3_5)),\qquad\text{ and} &(B(3_4,2,1),B(3_4,2,1)).
\end{align*}
In each case we will write $M=M_r=M_c$.  Because the row and column
solution sets are identical, some simplifications occur.

\subsubsection{$(B(4_3,3),B(4_3,3))$}
This Gram matrix corresponds to the conjectured maximal
determinant.  We will see that the decomposability analysis
completely determines the form of the determinant.

We may write $M=12I+B$ where $B$ is the rank-4 matrix
\begin{equation*}
B:=\begin{pmatrix}
3J_4 & -J_4 & -J_4 & -J_{4,3} \\
-J_4 & 3J_4 & -J_4 & -J_{4,3} \\
-J_4 & -J_4 & 3J_4 & -J_{4,3} \\
-J_{3,4} & -J_{3,4} & -J_{3,4} & 3J_3 \\
\end{pmatrix}.
\end{equation*}
If $x$ is a $\{-1,1\}$ row or column vector, define 
\begin{equation*}
a := \sum_{j=1}^4 x_j \qquad b := \sum_{j=5}^8 x_j \qquad
c := \sum_{j=9}^{12} x_j \qquad d := \sum_{j=13}^{15} x_j.
\end{equation*}
It follows that $a$, $b$, and $c$ are even and that $d$ is odd, and that
$-4\le a,b,c\le4$ and $-3\le d\le3$.
Equations~(\ref{rowquad}) and~(\ref{colquad}) both reduce to
\begin{equation*}
83=4(a^2+b^2+c^2+d^2)-(a+b+c+d)^2
\end{equation*}
for rows and columns 1 through 12.  It is a simple matter to find the
parity normalized solutions to this equation, subject to the above
constraints, using an exhaustive search.  The results are
\begin{equation*}
(a,b,c,d)=([-4, -4, 0], 1),\quad ([-4, -2, 2], -3),\quad ([-2, -2, 2], 3),
\text{ or } ([-2, 0, 4], -1)
\end{equation*}
where the notation $[a,b,c]$ indicates that any permutation of $a$, $b$, and
$c$ is allowed.  For rows and columns 13 through 15, the equation is
\begin{equation*}
75=4(a^2+b^2+c^2+d^2)-(a+b+c+d)^2,
\end{equation*}
and the solutions are
\begin{equation*}
(a,b,c,d)=(-4, -4, -4, 1),\quad ([-4, 0, 2], -1),\quad (-2, -2, -2, 3),
\text{ or } ([-2, 2, 4], 1).
\end{equation*}

Not all of these solutions for, say,  rows are compatible with each other.
Imposing the equation~(\ref{rowconsist}) for all inequivalent choices of $j$
and $k$ reduces the solution sets to
\begin{align*}
&\{(-4,-2,2,-3),(-2,0,4,-1)\} && \text{for rows 1--4,}\\
&\{(2,-4,-2,-3),(4,-2,0,-1)\} && \text{for rows 5--8,}\\
&\{(-2,2,-4,-3),(0,4,-2,-1)\} && \text{for rows 9--12,}\\
&\{(-4,-4,-4,1),(-2,-2,-2,3)\} && \text{for rows 13--15}.
\end{align*}
Note that the first three of these sets may be permuted.  Also the first
three components of all solution vectors may be simultaneously permuted.
Since such permutations give equivalent solutions, we
disregard them.  The column solutions must be also be compatible with each
other, and also with the row solutions, according to
equation~(\ref{rcconsist}).  The latter fixes a unique permutation of
the column solutions, namely
\begin{align*}
&\{(-4,2,-2,-3),(-2,4,0,-1)\} && \text{for columns 1--4,}\\
&\{(-2,-4,2,-3),(0,-2,4,-1)\} && \text{for columns 5--8,}\\
&\{(2,-2,-4,-3),(4,0,-2,-1)\} && \text{for columns 9--12,}\\
&\{(-4,-4,-4,1),(-2,-2,-2,3)\} && \text{for columns 13--15}.
\end{align*}

For each block of rows (or columns) two different row (or column) types
can occur.  By using~(\ref{decompPair}) we can determine precisely how many
times each type occurs in the block.  Consider rows 13 through 15.  If
any row has type $(-4,-4,-4,1)$, that row will have inner product 11 or
15 with any other row of the same type.   Hence a second row of type
$(-4,-4,-4,1)$ is not allowed as
the inner product must be 3.  On the other hand, a type $(-4,-4,-4,1)$ row
will have inner product
7 with any row of type $(-2,-2,-2,3)$, which is also disallowed.  Therefore
all three of the rows 13--15 must be of type $(-2,-2,-2,3)$.  Now the 
matrix formed by the last three rows
can be regarded as composed of three $3\times4$ blocks followed by a
$3\times3$ block
consisting entirely of 1s.  The $3\times4$ blocks must have row sums $-2$ and
column sums either $-3$ or $-1$.  The row-sum condition implies that there
is a single 1 in each row, while the column-sum condition implies that
there can be no more than a single 1 in each column.  Thus there is one
column with column sum $-3$ and three with column sum $-1$ in each block.
We may choose to place the $-3$ type column first in each block, and
likewise for rows.  Then the structure of the matrix is given by
\begin{align}\label{rowstruct}
&(-4,-2,2,-3) && \text{row 1}\notag\\
&(-2,0,4,-1) && \text{rows 2, 3, and 4}\notag\\
&(2,-4,-2,-3) && \text{row 5}\notag\\
&(4,-2,0,-1) && \text{rows 6, 7, and 8}\\
&(-2,2,-4,-3) && \text{row 9}\notag\\
&(0,4,-2,-1) && \text{rows 10, 11, and 12}\notag\\
&(-2,-2,-2,3) && \text{rows 13, 14, and 15}\notag,
\end{align}
and
\begin{align}\label{colstruct}
&(-4,2,-2,-3) && \text{column 1}\notag\\
&(-2,4,0,-1) && \text{columns 2, 3, and 4}\notag\\
&(-2,-4,2,-3) && \text{column 5}\notag\\
&(0,-2,4,-1) && \text{columns 6, 7, and 8}\\
&(2,-2,-4,-3) && \text{column 9}\notag\\
&(4,0,-2,-1) && \text{columns 10, 11, and 12}\notag\\
&(-2,-2,-2,3) && \text{columns 13, 14, and 15}\notag.
\end{align}

\subsubsection{Uniqueness of the decomposition}
The matrix found by Smith and Cohn exhibits the structure derived in the
previous section.  Here we will show that the matrix emerges quite simply
from the structure and that there is essentially only one such matrix.

The matrix must have the structure
\begin{equation}\label{matstruct}
R=\begin{pmatrix}
A_1 & B_1 & C & D_1 \\
C & A_2 & B_2 & D_2 \\
B_3 & C & A_3 & D_3 \\
E_1 & E_2 & E_3 & J_3 \\
\end{pmatrix},
\end{equation}
where $A_j$, $B_j$, and $C$ are $4\times4$ sub-matrices; $D_j$ are $4\times3$
sub-matrices; and $E_j$ are $3\times4$ sub-matrices.  The matrices $A_j$
have row and column sums $(-4,-2,-2,-2)$; $B_j$ have row and column sums
$(-2,0,0,0)$; and $C$ has row and column sums $(2,4,4,4)$ (which determines
$C$ completely and is why we need not distinguish different $C_j$).  The
matrices $D_j$ have row sums $(-3,-1,-1,-1)$ and column sums $(-2,-2,-2)$,
while for the $E_j$ the row sums are $(-2,-2,-2)$ and the column sums are
$(-3,-1,-1,-1)$.

The row and column sums imply that
\begin{equation*}
A_j\sim A:=
\begin{pmatrix}
- & - & - & - \\
- & 1 & - & - \\
- & - & 1 & - \\
- & - & - & 1 \\
\end{pmatrix},
\quad B_j\sim B:=
\begin{pmatrix}
1 & - & - & - \\
- & - & 1 & 1 \\
- & 1 & - & 1 \\
- & 1 & 1 & - \\
\end{pmatrix},
\quad C=
\begin{pmatrix}
- & 1 & 1 & 1 \\
1 & 1 & 1 & 1 \\
1 & 1 & 1 & 1 \\
1 & 1 & 1 & 1 \\
\end{pmatrix},
\end{equation*}
\begin{equation*}
D_j\sim D:=
\begin{pmatrix}
- & - & -  \\
1 & - & -  \\
- & 1 & -  \\
- & - & 1  \\
\end{pmatrix}, \qquad\text{ and }
\qquad E_j\sim E:=
\begin{pmatrix}
- & 1 & - & - \\
- & - & 1 & - \\
- & - & - & 1 \\
\end{pmatrix},
\end{equation*}
where ``$\sim$'' means equivalent up to permutation of the last three rows
or columns.  The structure~(\ref{rowstruct}), (\ref{colstruct}),
(\ref{matstruct}) is preserved under permutation of rows 2--4, of rows 6--8,
of rows 10--12, of rows 13--15, or of the corresponding columns.  We use
the freedom to permute columns to fix $A_1=A$, $B_1=B$, and $D_1=D$
in~(\ref{matstruct}).  Then use the freedom to permute rows to fix $A_2=A$,
$B_3=B$, and $E_1=E$.  Finally, another column permutation fixes $B_2=B$.
Thus rows 1--4 and columns 1--4 are completely constructed and have the
correct pairwise inner products.  In order for these to have the right
inner products with rows 5--8 and columns 5--8, we must have $D_2=D$ and
$E_2=E$.  For rows 13--15 to have the desired inner products with rows 6--8
necessitates $E_3=E$.  Similarly columns 13--15 have the correct inner products
with columns 2--4 only if $E_3=E$.  It then follows that $A_3=A$, and the
matrix $R$ such $RR^{\rm T}=R^{\rm T}R=M$ is completely determined up to
permutations consistent with the structure of $M$, and negation of the
entire matrix.  Consequently there is essentially only one D-optimal design
of order 15 (subject to the proof of non-decomposability of the other three
candidate Gram matrices, which follows).

\subsubsection{$(B(6,3,2_3),B(6,3,2_3))$}
The candidate matrix may be written $M=12I+B$ where $B$ is a rank-5 matrix.
Row and column vectors are subdivided in to 5 sub-vectors whose element sums
are denoted $a$, $b$, $c$, $d$, and $e$.  The sum $b$ is odd whereas the
other four are even.  The inequalities $-6\le a\le6$, $-3\le b\le3$, and
$-2\le c,d,e\le2$ must hold.

For rows or columns 1 through 6, equations~(\ref{rowquad}) and
(\ref{colquad}) reduce to
\begin{equation*}
99=4(a^2+b^2+c^2+d^2+e^2)-(a+b+c+d+e)^2.
\end{equation*}
Exhaustive search yields the solution set
\begin{align*}
(a,b,c,d,e)\in\{&(-6,-1,0,0,0),\ (-4,1,[0,2,2]),\ (-4,3,-2,-2,-2),\\
&(-2,-3,2,2,2),\ (2,-3,[-2,2,2]),\ (4,-3,0,0,0),\\
&(4,1,[-2,-2,0]),\ (6,-1,[0,2,2]),\ (6,3,0,0,0)\}.
\end{align*}
For rows or columns 7 through 9, equations~(\ref{rowquad}) and
(\ref{colquad}) reduce to
\begin{equation*}
75=4(a^2+b^2+c^2+d^2+e^2)-(a+b+c+d+e)^2
\end{equation*}
which has the solution set
\begin{align*}
(a,b,c,d,e)\in\{&(-6,-3,[-2,0,0]),\ (-6,1,-2,-2,-2),\ (-4,-1,[0,0,2]),\\
&(-2,3,[-2,-2,0]),\ (0,3,-2,-2,-2),\ (2,-3,2,2,2),\\
&(4,1,[-2,0,2])\}.
\end{align*}
For rows or columns 10 through 15, equations~(\ref{rowquad}) and
(\ref{colquad}) reduce to
\begin{equation*}
67=4(a^2+b^2+c^2+d^2+e^2)-(a+b+c+d+e)^2
\end{equation*}
which has the solution set
\begin{align*}
(a,b,c,d,e)\in\{&(-4,-1,[-2,-2,2]),\ (-2,1,[-2,2,2]),\ (0,-3,[0,2,2]),\\
&(2,-3,[0,0,2]),\ (2,1,[-2,-2,2]),\ (4,3,[-2,2,2]).
\end{align*}

If for a given row 1--6 solution the compatibility condition~(\ref{rowconsist})
fails to hold  for all row 7--9 solutions, we drop that solution.  Likewise
if a given row 7--9 solution is incompatible with all row 1--6 solutions,
we do not consider that solution any further.  The compatibility condition
in this case reduces to
\begin{equation*}
-33=x_7^{\rm T}(4I_5-J_5)x_1
\end{equation*}
where $x_j$ is a solution for row $j$ taken from the appropriate list above.
Failure to find any $x_7$ such that this equation holds when $x_1$ is
one of $(-4,1,[0,2,2])$, $(-4,3,-2,-2,-2)$, $(-2,-3,2,2,2)$, or
$(4,1,[-2,-2,0])$ eliminates those solutions.  Similarly the solutions
$x_7=(-2,3,[-2,-2,0])$ and $(4,1,[-2,0,2])$ are eliminated.

Imposing compatibility of the remaining row 1--6 solutions
with the row 10--15 solutions allows us to eliminate the solutions
$x_{10}=(0,-3,[0,2,2])$ and $(2,1,[-2,-2,2])$, and at the same time,
the solutions $x_1=(4,-3,0,0,0)$ and $(6,3,0,0,0)$.  As a consequence,
$x_7=(2,-3,2,2,2)$ is also eliminated, since the only $x_1$ solution
with which it was compatible was $(6,3,0,0,0)$.

Compatibility of the row 7--9 solutions with the row 10--15 solutions
then eliminates $x_{10}=(-2,1,[-2,2,2])$.

Finally, one may check that the compatibility condition between, say, row
10 and row 12, eliminates all the remaining row 10--15 solutions as
possibilities.  Therefore, there is no decomposition of
$(B(6,3,2,2,2),B(6,3,2,2,2))$.

\subsubsection{$(B(3,3,3,3,3),B(3,3,3,3,3))$}
The equations~(\ref{rowquad}) and~(\ref{colquad}) for the sums
of row and column sub-vectors have the same form for all 15 rows and
all 15 columns,
namely
\begin{equation*}
75=4(a^2+b^2+c^2+d^2+e^2)-(a+b+c+d+e)^2,
\end{equation*}
where $a$, $b$, $c$, $d$, and $e$ are odd integers in the range $[-3,3]$.
The parity normalized solutions to this equation are
$(a,b,c,d,e)=([-3,-3,1,1,1])$ and $(a,b,c,d,e)=([-3,-1,-1,-1,3])$.  Since
all solutions have at least one variable set to $-3$, we may use our freedom
to permute columns to fix $a=-3$ in the row 1 solution.  Likewise, we
fix $b=\pm3$ in row 1.  Now if $(a,b,c,d,e)=(-3,-3,1,1,1)$ in row 1, then
the consistency equation~(\ref{rowconsist}) for rows requires that the
same solution hold for rows 2 and 3 and that each of rows 4--15 be one
two types: $(a,b,c,d,e)=(-3,3,-1,-1,-1)$ or $(a,b,c,d,e)=(3,-3,-1,-1,-1)$.
On the other hand, if $(a,b,c,d,e)=(-3,3,-1,-1,-1)$ in row 1, then the
same solution must hold for rows 2 and 3, and the solution
$(a,b,c,d,e)=(-3,-3,1,1,1)$ must obtain in rows 4--15.

The combination of these two facts and permutation symmetry imposes the
requirement on the five sets
of rows, $\{1,2,3\}$, $\{4,5,6\}$, $\ldots$, $\{13,14,15\}$, that within
each set all three rows must be of the same type, but that outside the set,
all rows must be of the other type.  As no more than two row sets may
be constructed under these constraints, whereas we require five,
no decomposition exists.

\subsubsection{$(B(3,3,3,3,2,1),B(3,3,3,3,2,1))$}
The impossibility of decomposition follows by an argument that is similar to
but more involved than that of the previous case.  We need only look at
the solution set for rows 1--12.  Equation~(\ref{rowquad}) becomes
\begin{equation*}
75=4(a^2+b^2+c^2+d^2+e^2+f^2)-(a+b+c+d+e+f)^2
\end{equation*}
where $a$, $b$, $c$, and $d$ are odd integers in the range $[-3,3]$, $e$
is an even integer in the range $[-2,2]$, and $f\in\{-1,1\}$.  The
parity normalized solutions are in the set
\begin{align*}
\{&([-1,-1,3,3],2,-1),\ ([-3,1,1,3],2,1),\ ([1,1,3,3],-2,-1),\\
&([-1,1,3,3],-2,1),\ ([-3,-3,1,1],0,1),\ ([-3,-1,-1,3],0,-1)\}.
\end{align*}

We now use the compatibility conditions among rows~(\ref{rowconsist}) to
prove the impossibility of a decomposition.  We leave the arithmetical
checking of the consistency of various pairs of solutions to the reader.
A requirement of decomposability
is that we be able to assign consistent solutions to the four sets of
rows $\{1,2,3\}$, $\ldots$, $\{10,11,12\}$.  If the set $\{1,2,3\}$ has
a row of type $(-1,-1,3,3,2,-1)$ then rows 4--12 must be of
type $(1,1,[-3,3],2,1)$.  But two rows, in different row sets,
which both are of
the latter type are not compatible with each other, so it is not possible
to construct four sets when one of them contains a row of type
$([-1,-1,3,3],2,-1)$.

Similarly, if the set $\{1,2,3\}$ contains a row of type $(-3,1,1,3,2,1)$,
then rows 4--12 must all be of type $(1,3,3,1,-2,-1)$.  Again, two rows
of this type in different row sets are not compatible, so $(-3,1,1,3,2,1)$
has been eliminated.

If either of the two types $(-1,1,3,3,-2,1)$ and $(1,1,3,3,-2,-1)$ appears
in the set of rows $\{1,2,3\}$, then rows 4--12 are of a suitable form of type
$([-3,-3,1,1],0,1)$ or type $([-3,-1,-1,3],0,-1)$.  This means that at least
two of the three sets $\{4,5,6\}$, $\{7,8,9\}$, $\{10,11,12\}$
must both contain a row of the first
type, or both contain a row of the second type.  However neither of the latter
two types can appear in two different row sets.  Hence the solutions
$([-1,1,3,3],-2,1)$ and $([1,1,3,3],-2,-1)$ are ruled out.  But since, as we
just argued, neither of the two remaining solution types, $([-3,-3,1,1],0,1)$
and $([-3,-1,-1,3],0,-1)$ can appear in two different row sets,
constructing four row sets is impossible.  Therefore there is no decomposition.

\subsection{Other orders}
In orders 3, 7, and 9, decompositions of the the unique largest candidate
Gram matrices have been known for some time.

\subsubsection{$n=11$}
Ehlich proved the non-decomposability of the four candidates of largest
determinant, and showed that the three candidates of the next largest
determinant do decompose.  We have applied our method to these seven
matrices, and confirmed his results.  As the method does not differ in
any substantial way from that applied to the case $n=15$, we omit the
details.

\subsubsection{$n=17$}
In~\cite{MK} eight of the nine reported candidate Gram matrices were shown to
be non-decomposable, thus establishing the maximality of the determinant
associated with the ninth.  Since our list contains two additional candidates,
we must, in order to verify that result, prove them non-decomposable as well.
Furthermore, one of the additional candidates has the same characteristic
equation as one of the original nine, so this pair must also be checked.

The three decompositions to be tested are then
\begin{align*}
&(S(5_2,-3_2,1_{12}),S(5_2,-3_2,1_{12})),\\
&(D(1;1,-3_6,1_8;3,1_{14}),D(1;1,-3_6,1_8;3,1_{14})),\\
&(S(5_2,-3_2,1_{12}),S(-3_8,1_8)).
\end{align*}
For the former two, follow the model of proof in~\cite{MK}.  (We have
checked these, as well as the original nine candidates.)  As the last
pair gives us a chance to illustrate the method when $M_r\ne M_c$, we
present the proof here.

We let $M_r=S(5_2,-3_2,1_{12})$ and $M_c=S(-3_8,1_8)$.  We compute the
sets of solutions to the quadratic equation~(\ref{rowquad}) for row 1
and for rows 2--3, and then show that the sets are incompatible.  The
equation for row 1 is
\begin{equation*}
97=a^2+b^2+c^2-6ab+2ac+2bc
\end{equation*}
where $a\in\{-1,1\}$ and $b$ and $c$ are even integers in the range $[-8,8]$.
The parity normalized solutions are $(a,b,c)=(1,-2,-8)$ and $(a,b,c)=(-1,6,2)$.
For rows 2 and 3, the equation is
\begin{equation*}
57=a^2+b^2+c^2-6ab+2ac+2bc
\end{equation*}
and the parity normalized solutions are $(-1,4,-8)$, $(-1,6,-2)$,
$(1,-6,8)$, $(1,-4,-2)$, and $(1,8,2)$.  The compatibility
constraint~(\ref{rowconsist})
\begin{equation*}
m_1\cdot m_2 - 16m_{1,2}=x_1^{\rm T}(M_c-16I)x_2
\end{equation*}
which reduces to
\begin{equation*}
101=a_1 a_2+b_1 b_2+c_1 c_2-3(a_1 b_2+b_1 a_2)+a_1 c_2+c_1 a_2+b_1 c_2+c_1 b_2
\end{equation*}
where $x_1^{\rm T}$ is represented by $(a_1,b_1,c_1)$ and $x_2^{\rm T}$
by $(a_2,b_2,c_2)$.  Taking all pairs of solutions from the sets
enumerated above, we find that none satisfy the constraint.  Therefore
no decomposition exists.

\subsubsection{$n=19$}
There are presently no known candidate Gram matrices.

\subsubsection{$n=21$}
As in the case $n=17$, in order to verify the result of~\cite{CKM} we
have three additional pairs of Gram matrices
to prove non-decomposable,
\begin{align*}
&(S(5,-3_5,1_{14}),S(5,-3_5,1_{14}),\\
&(S(-7,-3_2,1_{17}),S(-7,-3_2,1_{17})),\text{ and }\\
&(S(5,-3_5,1_{14}),S(-7,-3_2,1_{17})).
\end{align*}
The proof parallels that of the $n=17$ case, so we omit it here.

\subsubsection{$n=29$, $33$, and $37$}
All of the candidate Gram matrices enumerated in Table~\ref{table:29-33-37}
have been shown not to decompose.  As there are a great many matrices
to check, an automated procedure was used for this purpose which combines
some, but not all, aspects of the method illustrated in the case of $n=15$
with a
backtracking search.  (In fact, the result for $n=15$ was originally
obtained by the use of this algorithm. As its performance is far from
optimal, we will not describe it in detail here.  We hope to present an
improved version in a forthcoming publication.)

Combining the above with the results of~\cite{OSDS} and~\cite{OS} we are able
to place new bounds on $\mdn(n)$ for these three orders.  For $n=29$
we conclude that $320\times7^{12}\le\mdn(29)<342\times7^{12}$ or in
terms of Barba's bound~(\ref{Bone}),
$0.865001\ldots\le\md(29)/B(29)<0.92447$.  For $n=33$ we have
$441\times8^{14}\le\mdn(33)<485\times8^{14}$ or
$0.854667\ldots\le\md(33)/B(33)<0.939951$.  Finally, for $n=37$
we find $72\times9^{17}\le\mdn(37)<659\times9^{16}$ or
$0.936329\ldots\le\md(37)/B(37)<0.952224$.

\section{The range of the determinant function for $n=9$ and $n=11$}
Another application of these methods is the determination of the complete
range of the determinant function for $\{-1,1\}$-matrices.  By
negation of appropriate rows and columns, we can make the first
row of a matrix consist entirely
of 1s and the first column, except for the $(1,1)$ entry consist entirely
of $-1$s.  Such matrices are called {\em normalized.}  Adding
the first row of a normalized matrix to every other row
does not change the determinant, and produces a matrix whose first column
has a single 1 in the $(1,1)$ position and 0s elsewhere.  All entries in
rows 2 through $n$, where $n$ is the size of the matrix, are either 0 or 2.
Expansion of the determinant by minors on the first column then shows that
the determinant is equal to that of a $\{0,2\}$-matrix of size $n-1$.  Thus
an $n\times n$ $\{-1,1\}$-determinant is always divisible by $2^{n-1}$.
Furthermore, we have shown the problem of maximizing the determinant of
$n\times n$ $\{-1,1\}$-matrices to be equivalent to that of maximizing 
$(n-1)\times(n-1)$ $\{0,1\}$-matrices.

The problem of the range of the determinant function was studied by
Craigen~\cite{Cr} who asked for the complete list of integers, $d$, such
that $d$ is the determinant of some $(n-1)\times(n-1)$ $\{0,1\}$-matrix, or
equivalently, $2^{n-1}d$ is the determinant of some $n\times n$
$\{-1,1\}$-matrix.  We are focusing here on $\{-1,1\}$-matrices,
but it is convenient
to omit the factor $2^{n-1}$ in the discussion, so we will be careful always
to refer to {\em normalized} determinants (see Section~\ref{sect:intro}).

We ignore the sign of the determinant.  For $n\le7$ the
range of the normalized determinant function is $[0,\mdn(n)]$,
where the notation $[a,b]$ means all integer values between $a$ and $b$.
Craigen proved that in $n=8$, a gap appears for the first time.
In particular he showed that there are no normalized determinants in
the range $[28,31]$ whereas the order 8 Hadamard matrix has normalized
determinant 32.

Craigen also provided a table of normalized determinant values known to him
at the time.  For $n=8$ and $n=9$ we reproduce his lists in
Table~\ref{table:range}.  (For $n=8$ Craigen also lists $19$, which is almost
surely a misprint.)  For orders 10 and 11, Craigen gave partial lists.
In the table, we give what we believe to be complete lists of normalized
determinants.
\begin{table}
\begin{tabular}{llll}
$n=8$ & $n=9$ & $n=10$ & $n=11$\\
\hline
[0, 18] & [0, 40] & [0, 102] & [0, 268]\\
20 & 42 &[104, 105] & [270, 276]\\
24 & [44, 45] &108 & [278, 280]\\
32& 48 &110 & [282, 286]\\
& 56 &112 & 288\\
&& [116, 117] & 291\\
&& 120 & [294, 297]\\
&& 125 & 304\\
&& 128 & 312\\
&& 144 & 315\\
&&& 320
\end{tabular}
\caption{The range of the normalized determinant function for orders 8, 9,
10, and 11.\label{table:range}}
\end{table}

We have been able to exhibit a determinant for each of the values in the
table.  The method by which this was done is beyond the scope of this paper,
and we intend to describe it elsewhere.  What is of relevance to us here
is that for the two
odd orders, $n=9$ and $n=11$, we are able to prove that the lists are
complete, that is, there are no normalized determinant values other than
those given in the table.

The method is as follows: we set the threshold in our candidate Gram matrix
finding routine to the lowest gap in the range.  For $n=9$ this is
$41\times2^8$ and for $n=11$ it is $269\times2^{10}$.  The program then
produces a complete list of candidate Gram matrices corresponding to
determinants equaling or exceeding the threshold.  If none of these
determinants correspond to any of the gaps, then we have proved the claim.
On the other hand, if any of the candidate Gram matrices do correspond to
gaps, we must then show that these matrices do not decompose.

In the case $n=9$ we find candidate Gram matrices corresponding to the
following normalized determinants: 56~(1 matrix), 48~(4 matrices),
45~(1 matrix), 44~(2 matrices), and 42 (1 matrix).  Since none of these
correspond to gaps, completeness is proved.

\begin{table}
\begin{tabular}{l|cccccccccc}
det. & 270 & 271 & 272 & 273 & 274 & 275 & 276 & 278 & 279 & 280\\
mult. & 47 & 1 & 37 & 8 & 2 & 8 & 5 & 1 & 7 & 10\\
\hline
det. & 282 & 283 & 284 & 285 & 286 & 288 & 291 & 294 & 295 & 296\\
mult. & 2 & 1 & 2 & 2 & 1 & 41 & 1 & 2 & 1 & 1 \\
\hline
det. & 297 & 300 & 304 & 306 & 312 & 315 & 320 & 324\\
mult. & 1 & 3 & 2 & 1 & 1 & 1 & 3 & 4\\
\end{tabular}
\caption{Multiplicities of the candidate Gram matrices corresponding to
the listed normalized determinants of order 11.  The threshold used is
$269$.\label{table:11dets}}
\end{table}
For $n=11$ we find matrices corresponding to the determinants listed in
Table~\ref{table:11dets}.  There are 196 matrices in total.  All of the
determinant values correspond to gaps except for
three of them, namely 300, 306, and 324.  The non-decomposability of the
four matrices with determinant $(324\times2^{10})^2$ was discussed
previously.  To establish the desired result, one need only show that the
three matrices with determinant $(300\times2^{10})^2$ and the matrix
with determinant $(306\times2^{10})^2$ do not decompose.  These matrices
are
\begin{equation*}
\begin{pmatrix}
11&3&3&-&-&&&\\
3&11&-&3&-&&&\\
3&-&11&-&3&&-J_{5,6}&\\
-&3&-&11&3&&&\\
-&-&3&3&11&&&\\
&&&&&&&\\
&&\mspace{-18.0mu}-J_{6,5}\mspace{-18.0mu}&&&&B(3_2)&\\
&&&&&&&\\
\end{pmatrix},
\begin{pmatrix}
11&3&3&-&&&\\
3&11&-&3&&&\\
3&-&11&3&&-J_{4,7}&\\
-&3&3&11&&&\\
&&&&&&\\
&&\mspace{-18.0mu}-J_{7,4}\mspace{-18.0mu}&&&B(3_2,1)&\\
&&&&&&&\\
\end{pmatrix},
\end{equation*}
\begin{equation*}
\begin{pmatrix}
11&3&3&-&-&&&\\
3&11&-&3&-&&&\\
3&-&11&-&3&&-J_{5,6}&\\
-&3&-&11&3&&&\\
-&-&3&3&11&&&\\
&&&&&&&\\
&&\mspace{-18.0mu}-J_{6,5}\mspace{-18.0mu}&&&&B(4,1_2)&\\
&&&&&&&\\
\end{pmatrix},\text{ and }
\begin{pmatrix}
11&3&3&3&-&&&\\
3&11&-&-&3&&&\\
3&-&11&-&-&&-J_{5,6}&\\
3&-&-&11&-&&&\\
-&3&-&-&11&&&\\
&&&&&&&\\
&&\mspace{-18.0mu}-J_{6,5}\mspace{-18.0mu}&&&&B(3,1_3)&\\
&&&&&&&\\
\end{pmatrix}.
\end{equation*}
There are no unusual features in the proofs of non-decomposability,
so we omit them.

\section{Conclusion and outlook}
We have proved that the maximal determinant of a $15\times15$ $\{-1,1\}$-matrix
is $25515\times2^{14}$ and confirmed the known maximal values for the other
odd orders which do not attain the Barba bound, up to order $21$, except for
$19$.  We have also established new upper bounds on the maximal determinant
in orders $29$, $33$, and $37$.

Convincing conjectures exist for the maximal determinants for orders
$n=19$~\cite{Sm}, $22$, $23$~\cite{OSDS}, and $37$~\cite{OS}.  We hope that
the methods of this paper can be extended to handle at least some of
these cases.

Order $19$ may be tractable using the current method with
efficiency improvements in the computer code, and perhaps parallelization.
For higher orders, $n\equiv3\bmod 4$, there are excellent lower bounds on
the maximal determinant, but these can possibly be improved.  It
is unlikely that our method can be applied to these cases without a new
algorithm.  

The smallest even order for which the maximal determinant is
unknown is $22$, which may well be tractable, although it is too early to
tell until experiments are done.  

The biggest hope is for orders
$n\equiv1\bmod4$, the lowest open case of which is $n=29$.  The major hurdle
to be overcome here is the inefficiency of the program that searches for
candidate Gram matrices.  Efforts are underway to try to improve this.
In general, the $n\equiv1\bmod 4$ cases
are the most amenable to our method, and the major determinant of success
is how close the actual maximum lies to the Barba bound.  Thus we are
especially optimistic about the case $n=37$.

\section*{Acknowledgments}
I carried out much of this work as a Research Associate in
the Department of Mathematics and Statistics at the University of Melbourne,
and gratefully acknowledge the support of Tony Guttmann and of the Australian
Research Council.  A stimulating collaboration with
Judy-anne Osborn extending this work has produced some
new insight into the original result.  In particular I thank her for
a discussion in which a key property of the lexicographic ordering used
in the Gram matrix finding algorithm was understood.
I used \emph{Mathematica} extensively during this
project and also thank Indiana University for the use of its Sun E10000
and IBM RS/6000 SP computing platforms.

\end{document}